# SOME REMARKS ON METRICS INDUCED BY A FUZZY METRIC

R. ROOPKUMAR AND R. VEMBU

ABSTRACT. We introduce a crisp metric $d_M$ as the common limit of two different nets $(\Delta_{M,\lambda})$ and $(\delta_{M,\lambda})$ of crisp metrics induced by a fuzzy metric $M$ and prove that the existence of each of these limits is equivalent to that of the other and it is characterized by another condition on the original fuzzy metric $M$. We also derive some of the properties of these approximate metrics $\Delta_\lambda$ and $\delta_\lambda$. On the other hand, for a given a crisp metric $d$, establish that the fuzzy metric representing $M_d$ with values in $\{0,1\}$ and $d$ are compatible with the same topology. Further, we prove that if a crisp metric $d$ induces a fuzzy metric $M_d$, then all the approximate crisp metrics $\Delta_{M,\lambda}$ and $\delta_{M,\lambda}$ induced by this fuzzy metric are equal to the original metric $d$.

## 1. INTRODUCTION

The concept of fuzzy metric was introduced by Karamosil and Michálek [10], with motivations from the notions of statistical metrics [11, 21, 12, 20] and probabilistic metrics [13, 15]. Then, different notions of fuzzy metrics have been introduced in the literature. To mention a few, pseudo-quasi metric (respectively, fuzzy pseudo metric) which assigns distance between two fuzzy subsets (respectively, fuzzy points) as a nonnegative real number in [5, 4]; the fuzzy metric in [9] assigns the distance between two points in a set as a fuzzy real number; George and Veeramani [6] slightly altered the definition of fuzzy metric introduced in [10] to achieve some topological properties of the fuzzy metric space; the concept of metrics of fuzzy subsets is discussed in [2]; Das [3] introduced the notion of $L$-fuzzy numbers and defined distance between two points as an $L$-fuzzy number; Savchenko and Zarichnyi [19] modified the axiom of triangle inequality and called the new fuzzy metric as ultra fuzzy metric; and in [22] a fuzzy distance between two points is defined as an $L$-fuzzy number. A comparison of different fuzzy metrics can be seen from [18]. Among the different notions of fuzzy metrics in the literature, the fuzzy metric introduced by Karamosil, Michálek [10] and modified by George and Veeramani [6] attracted more researchers to pursue their research in this area of research.

When there is an uncertainty in measuring the distance between two points, the concept of fuzzy metric is involved. In view of practical purpose, it is necessary to find a crisp metric from the fuzzy metric, in a convincing manner. This is the motivation for introducing crispification of fuzzy metric. In pure mathematical point of view, if a fuzzy metric is modified as a crisp metric, then it is natural to ask a procedure to get back to the original fuzzy metric. So, we also discuss the fuzzification of crisp metric. Already, there are uncountable number of fuzzy metrics induced by a crisp metric in the literature [6]. However, none of them is a function with range $\{0,1\}$, which is in contrast with the identification of a crisp subset with a fuzzy subset. So, we introduce a fuzzy metric with range $\{0,1\}$, representing the given crisp metric. Finally, we prove that these two processes are consistent in the sense that for a crisp metric $d$, all the approximate metrics obtained from the fuzzy metric with range $\{0,1\}$ induced by $d$ are equal to the original metric.

## 2. Preliminaries

**Definition 2.1.** [10, Definition 7] *A fuzzy metric $M$ on a set $X$ is a fuzzy set of $X \times X \times \mathbb{R}$ satisfying the following conditions.*

(KM1) $M(x, y, t) = 0$ *for all $x, y \in X$ and all $t \leq 0$,*
(KM2) $M(x, y, t) = 1$ *for all $t > 0$ if and only if $x = y$,*
(KM3) $M(x, y, t) = M(y, x, t)$ *for all $x, y \in X$ and all $t \in \mathbb{R}$,*
(KM4) $M(x, z, t + s) \geq S(M(x, y, t), M(y, z, s))$, *where $S$ is a measurable binary real function defined on $[0,1] \times [0,1]$ taking its value on $[0,1]$ and such that $S(1,1) = 1$,*
(KM5) *For every pair $(x, y) \in X \times X$, $M(x, y, t)$ is a left continuous and non-decreasing function of $t$ such that $\lim_{t \to \infty} M(x, y, t) = 1$.*

(1)      Here, $M(x, y, t)$ is interpreting the degree of the truth value of the statement "*the distance between $x$ and $y$ is less than $t$*".

George and Veeramani[6] altered the above definition by changing the domain $X \times X \times \mathbb{R}$ and codomain $\mathbb{R}$ of the fuzzy metric, respectively, by $X \times X \times (0, \infty)$ and $(0,1]$, the measurable binary real function $S$ by a $t$-norm, left continuity of $M(x, y, \cdot)$ by continuity and by ignoring the monotonicity of $M(x, y, \cdot)$ and limit of $M(x, y, \cdot)$ at infinity. It is interesting to see that the monotonicity of $M(x, y, \cdot)$ is proved in [7] from the definition of fuzzy metric defined as follows.

**Definition 2.2.** [6] *The $3$-tuple $(X, M, *)$ is said to be a fuzzy metric space if $X$ is an arbitrary set, $*$ is a continuous $t$-norm and $M$ is a fuzzy set on $X \times X \times (0, \infty)$ satisfying the following conditions:*

(GV1) $M(x, y, t) > 0$,
(GV2) $M(x, y, t) = 1$ *for all $t$ if and only if $x = y$,*
(GV3) $M(x, y, t) = M(y, x, t)$,
(GV4) $M(x, y, t) * M(y, z, s) \leq M(x, z, t + s)$,
(GV5) $M(x, y, \cdot) : (0, \infty) \to [0, 1]$ *is continuous,*

$x, y, z \in X$ *and $t, s > 0$, where a continuous $t$-norm is a binary operation $*$ defined on $[0, 1]$ such that $([0, 1], *)$ is a topological monoid with unit $1$ and $a * b \leq c * d$ whenever $a \leq c$ and $b \leq d$, $a, b, c, d \in [0, 1]$.*

Although, the notion of $t$-norms is generalizing *minimum*, throughout this paper we prefer to use minimum as the $t$-norm, and hence we denote a fuzzy metric space simply by $(X, M)$.

In this paper, for our discussion we use another simple modification of the fuzzy metric, which is obtained by including one more natural condition to the set of conditions given Definition 2.1, as follows.

**Definition 2.3.** *A fuzzy metric $M$ on an arbitrary set $X$ is a function $M : X \times X \times \mathbb{R} \to [0,1]$ satisfying the following conditions.*

(KM1) $M(x, y, t) = 0$ *for all $x, y \in X$ and all $t \leq 0$,*
(KM2) $M(x, y, t) = 1$ *for all $t > 0$ if and only if $x = y$,*
(KM3) $M(x, y, t) = M(y, x, t)$ *for all $x, y \in X$ and all $t \in \mathbb{R}$,*

(KM4) $M(x, z, t + s) \geq \min\{M(x, y, t), M(y, z, s)\}$,
(KM5) *If $x, y \in X$ and $x \neq y$, then $M(x, y, t)$ is a left continuous function of $t$ such that*
$$\lim_{t \to \infty} M(x, y, t) = 1,$$
(SDP) *If $x \neq y$, then* $\lim_{t \to 0} M(x, y, t) = 0$.

Since the distance between $x$ and $y$ must be positive whenever $x \neq y$, as the distance parameter $t$ tends to 0, the truth value of the statement "*the distance between $x$ and $y$ is less than $t$*" is expected to tend to zero. This property of separating distinct points is included as the condition (SDP).

**Example 2.4.** *For any set $X$ define $M : X \times X \times \mathbb{R} \to [0, 1]$ by*

$$M(x, y, t) = \begin{cases} 0 & \text{if} \quad t \leq 0, \text{ for all } x, y \in X, \\ 1 & \text{if} \quad x = y \text{ and } t > 0, \\ t & \text{if} \quad x \neq y \text{ and } t \in (0, \frac{1}{4}] \cup (\frac{3}{4}, 1], \\ \frac{1}{2} & \text{if} \quad x \neq y \text{ and } t \in (\frac{1}{4}, \frac{3}{4}], \\ 1 & \text{if} \quad t > 1, \text{ for all } x, y \in X. \end{cases}$$

*Then $M$ is a fuzzy metric.*

For an arbitrarily fixed distinct pair $x, y \in X$, the graph of $M(x, y, t)$ is the following.

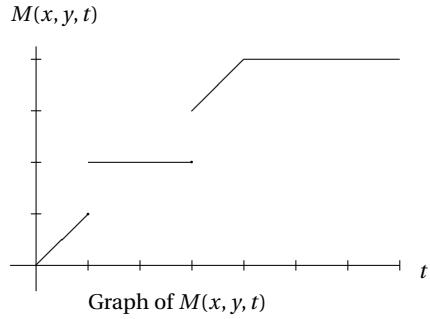

Graph of $M(x, y, t)$

**Example 2.5.** *For any set $X$ define $M : X \times X \times \mathbb{R} \to [0, 1]$ by*

$$M(x, y, t) = \begin{cases} 0 & \text{if} \quad t \leq 0, \text{ for all } x, y \in X, \\ 1 & \text{if} \quad x = y \text{ and } t > 0, \\ \frac{1}{2} & \text{if} \quad x \neq y \text{ and } 0 < t \leq \frac{1}{2}, \\ 1 & \text{if} \quad t > \frac{1}{2}, \text{ for all } x, y \in X. \end{cases}$$

*Then $M$ is a fuzzy metric according to Definition 2.1 whereas it is not a fuzzy metric according to Definition 2.3 as it does not satisfy Condition (SDP).*

In [6], open balls in a fuzzy metric space are defined and a topology is generated by the collection of all open balls as a basis.

**Definition 2.6.** [6, Definition 3.1] *Let $(X, M)$ be a fuzzy metric space. The open ball $B(x, r, t)$ for $t > 0$ with centre $x \in X$ and radius $r$, $0 < r < 1$ is defined as $B(x, r, t) = \{y \in X : M(x, y, t) > 1 - r\}$.*

It is unusual to restrict the radius of an open ball in $0 < r < 1$. Since the distance is represented by the parameter $t$, and the degree of the fact that the distance between $x$ and $y$ is less than $t$ is $M(x, y, t)$, it is better to rewrite the above definition as follows.

**Definition 2.7.** *Let* $(X, M)$ *be a fuzzy metric space. We define open ball* $B(x, r, \epsilon)$ *with centre* $x \in X$, *radius* $r > 0$ *and a parameter of fuzziness* $0 < \epsilon < 1$ *as* $B(x, r, \epsilon) = \{y \in X : M(x, y, r) > 1 - \epsilon\}$.

So, for every $y \in B(x, r, \epsilon)$, the truth value of the statement that *the distance between $x$ and $y$ is less than $r$* is $1 - \epsilon$. This change is only on the interpretation of the open ball $B(x, r, \epsilon)$, and it does not make any difference between the two collections $\{B(x, r, t) : 0 < r < 1, t > 0\}$ and $\{B(x, r, \epsilon) : r > 0, 0 < \epsilon < 1\}$.

## 3. Crispification of a fuzzy metric

First, we introduce upper $\lambda$-metrics and lower $\lambda$-metrics by which the distance between two points are calculated upto the degree of *correctness parameter* $\lambda$.

**Definition 3.1** ($\lambda$-metric induced by a fuzzy metric)**.** *Let* $(X, M)$ *be a fuzzy metric space and* $\lambda \in (0, 1)$. *Let* $\Delta_{M;\lambda}$ *and* $\delta_{M;\lambda}$ *be defined by*

$$\Delta_{M;\lambda}(x, y) = \inf\{t \in \mathbb{R} : M(x, y, t) > \lambda\}$$
$$\delta_{M;\lambda}(x, y) = \sup\{t \in \mathbb{R} : M(x, y, t) < \lambda\},$$

$\forall (x, y) \in X \times X$.

**Theorem 3.2.** *Let* $(X, M)$ *be a fuzzy metric space. For each* $\lambda \in (0, 1)$, $\delta_{M;\lambda}$ *and* $\Delta_{M;\lambda}$ *are metrics on* $X$.

*Proof.* Since, for each pair $x, y$, $M(x, y, t)$ is an increasing function of $t$, we observe that $\{t \in \mathbb{R} : M(x, y, t) > \lambda\}$ is an interval with left end $\Delta_{M;\lambda}(x, y)$ and right end $+\infty$.

Clearly, $\Delta_{M;\lambda}(x, y) \geq 0$ and $\Delta_{M;\lambda}(x, y) = \Delta_{M;\lambda}(y, x)$ for all $x, y \in X$. If $x = y$, then $M(x, y, t) = 1$ for all $t > 0$, which implies that $\{t : M(x, y, t) > \lambda\} = (0, \infty)$ and hence $\Delta_{M;\lambda}(x, y) = 0$. Conversely, suppose that $\Delta_{M;\lambda}(x, y) = 0$ and $x \neq y$. Since $M(x, y, t)$ is right continuous at 0 (by condition $(SDP)$ and $M(x, y, 0) = 0$, there exists $t_0 > 0$ such that $M(x, y, t_0) < \lambda$; this implies that $t_0 \notin \{t : M(x, y, t) > \lambda\}$ and hence $\Delta_{M;\lambda}(x, y) \geq t_0 > 0$, which is a contradiction. Thus, we have proved that $\Delta_{M;\lambda}(x, y) = 0$ if and only if $x = y$.

Let $x, y, z \in X$. If any two of $x, y$ and $z$ are equal, then it follows that $\Delta_{M;\lambda}(x, z) \leq \Delta_{M;\lambda}(x, y) + \Delta_{M;\lambda}(y, z)$. So we assume that $x, y$ and $z$ are pairwise distinct. By the observation in the first line of this proof, we have

$$M\left(x, y, \Delta_{M;\lambda}(x, y) + \frac{\epsilon}{2}\right) > \lambda, \ M\left(y, z, \Delta_{M;\lambda}(y, z) + \frac{\epsilon}{2}\right) > \lambda;$$

and hence $M(x, z, \Delta_{M;\lambda}(x, y) + \Delta_{M;\lambda}(y, z) + \epsilon) > \lambda$. This implies that

$$\Delta_{M;\lambda}(x, y) + \Delta_{M;\lambda}(y, z) + \epsilon \in \{t : M(x, z, t) > \lambda\}$$

which shows that $\Delta_{M;\lambda}(x, z) \leq \Delta_{M;\lambda}(x, y) + \Delta_{M;\lambda}(y, z) + \epsilon$. As this is true for all $\epsilon > 0$ we have $\Delta_{M;\lambda}(x, z) \leq \Delta_{M;\lambda}(x, y) + \Delta_{M;\lambda}(y, z)$ and hence $\Delta_{M;\lambda}$ is a metric.

By employing a similar set of arguments we can show that $\delta_{M;\lambda}$ is a metric. □

Hereafter, we call $\Delta_{M;\lambda}$ and $\delta_{M;\lambda}$, respectively, as the upper $\lambda$-metric and the lower $\lambda$-metric induced by the given fuzzy metric $M$. It is easy to observe that for each $\lambda \in (0, 1)$, we have $\delta_{M;\lambda} \leq \Delta_{M;\lambda}$. In fact, for the fuzzy metric $M$ given in Example 2.5, we can see that $\delta_{M;\frac{1}{2}}(x, y) = \frac{1}{4} < \frac{3}{4} = \Delta_{M;\frac{1}{2}}(x, y)$, for an arbitrary $x, y \in X$ with $x \neq y$. The following theorem provides a necessary and sufficient condition for getting equality between $\delta_{M;\lambda}$ and $\Delta_{M;\lambda}$

**Theorem 3.3.** *Let $M$ be a fuzzy metric on $X$. Then, for any $x, y \in X$, $\delta_{M;\lambda}(x, y) = \Delta_{M;\lambda}(x, y)$ if and only if the set $\{t : M(x, y, t) = \lambda\}$ contains at most one element.*

*Proof.* Let $S = \{t : M(x,y,t) = \lambda\}$, $S_l = \{t : M(x,y,t) < \lambda\}$ and $S_u = \{t : M(x,y,t) > \lambda\}$. Suppose that $S$ contains two elements $t_1$ and $t_2$ with $t_1 < t_2$, then $t_1 \notin S_l$, $t_2 \notin S_u$ and hence $\delta_{M;\lambda}(x,y) \leq t_1 < t_2 \leq \Delta_{M;\lambda}(x,y)$, which implies that $\delta_{M;\lambda}(x,y) \neq \Delta_{M;\lambda}(x,y)$.

Conversely, if $\delta_{M;\lambda}(x,y) < \Delta_{M;\lambda}(x,y)$, then we choose a real number $s$ such that $\delta_{M;\lambda}(x,y) < s < \Delta_{M;\lambda}(x,y)$. Therefore, it follows that $s \notin S_l \subseteq (-\infty, \delta_{M;\lambda}(x,y)]$ and $s \notin S_u = [\Delta_{M;\lambda}(x,y), +\infty)$ and hence $s \in S$. Thus $S$ contains uncountable elements, which completes the proof. □

**Lemma 3.4.** *Let $(X,M)$ be a fuzzy metric space. If $0 < \lambda_1 < \lambda_2 < 1$, then $\delta_{M;\lambda_1} \leq \delta_{M;\lambda_2}$ and $\Delta_{M;\lambda_1} \leq \Delta_{M;\lambda_2}$.*

*Proof.* Since $\lambda_1 < \lambda_2$, we have $\{t : M(x,y,t) < \lambda_1\} \subseteq \{t : M(x,y,t) < \lambda_2\}$ and $\{t : M(x,y,t) > \lambda_1\} \supseteq \{t : M(x,y,t) > \lambda_2\}$. So, it follows that $\delta_{M;\lambda_1}(x,y) \leq \delta_{M;\lambda_2}(x,y)$ and $\Delta_{M;\lambda_1}(x,y) \leq \Delta_{M;\lambda_2}(x,y)$. □

**Remark 3.5.** *In general, $\Delta_{M;\lambda_1}$ and $\Delta_{M;\lambda_2}$ need not be the same for distinct $\lambda_1$ and $\lambda_2$. Moreover, they need not be compatible with a same topology. Similar, statement is true for $\delta_{M;\lambda_1}$ and $\delta_{M;\lambda_2}$.*

**Example 3.6.** *Let $X = [0,1]$. For $x,y \in X$ and $t \in \mathbb{R}$, if*

$$M(x,y,t) = \begin{cases} 0 & \text{if } -\infty < t \leq 0 \\ 1 & \text{if } x = y \text{ and } t > 0 \\ \frac{3t}{4(t+|x-y|)} & \text{if } x \neq y \text{ and } 0 < t \leq 2 \\ 1 & \text{if } x \neq y \text{ and } 2 < t, \end{cases}$$

*then $(X,M)$ is a fuzzy metric.*

Let $\lambda_1 = \frac{1}{2}$ and $\lambda_2 = \frac{7}{8}$. For $x \neq y$, we have

$$\delta_{M;\lambda_1}(x,y) = \sup\left\{t : \frac{3t}{4(t+|x-y|)} < \frac{1}{2}\right\} = \sup\{t : t < 2|x-y|\} = 2|x-y|,$$

and $\Delta_{M;\lambda_1}(x,y) = \inf\{t : t > 2|x-y|\} = 2|x-y|$. Thus, the topology compatible with $\Delta_{M;\lambda_1} = \delta_{M;\lambda_1}$ is the relative topology on $[0,1]$ inherited from the usual topology on $\mathbb{R}$.

But, for $x \neq y$, $\Delta_{M;\lambda_2}(x,y) = \sup\left\{t : M(x,y,t) < \frac{7}{8}\right\} = \sup(-\infty, 2] = 2$ and $\delta_{M;\lambda_2}(x,y) = \inf\left\{t : M(x,y,t) > \frac{7}{8}\right\} = \inf(2,\infty) = 2$, and hence the topology compatible with $\Delta_{M;\lambda_2} = \delta_{M;\lambda_2}$ is the discrete topology.

Our next questions are the following.
(1) Do $\lim_{\lambda \to 1} \Delta_{M;\lambda}$ and $\lim_{\lambda \to 1} \delta_{M;\lambda}$ exist?
(2) If so, does the existence of one limit imply that of the other?
(3) If the limits exist, are they equal?
(4) If the limits exist, are they metrics?

If $M$ is induced by a metric $d$ on $X$ (as in Theorem 4.2), then the limit exists and it is equal to $d$. But in general, the limit does not exist.

**Example 3.7.** *For $x,y \in \mathbb{R}$, if $M(x,y,t) = \begin{cases} \frac{t}{t+|x-y|} & t > 0 \\ 0 & t \leq 0 \end{cases}$, then $M$ is a fuzzy metric on $\mathbb{R}$. Let us take $x = 2$ and $y = 1$. Then $M(x,y,t) = \frac{t}{t+1}$ and $\Delta_{M;\lambda}(x,y) = \frac{\lambda}{1-\lambda}$. Thus $\Delta_{M;\lambda}(x,y) \to \infty$ as $\lambda \to 1$.*

To characterize the fuzzy metrics, for which the limits $\lim_{\lambda \to 1} \Delta_{M;\lambda}$ and $\lim_{\lambda \to 1} \delta_{M;\lambda}$ exist, we introduce a particular class of fuzzy metric spaces, which satisfy the *finite distance* condition

(*FD*) *For every pair* $(x, y)$, *there exists* $t_{x,y}$ *such that* $M(x, y, t_{x,y}) = 1$,

which is reflecting the natural expectation that distance between any two points is finite and it was introduced in Page 342 of [10].

**Definition 3.8** (*FD*-fuzzy metric space). *A fuzzy metric space* $(X, M)$ *(as in Definition 2.3) is said to be an FD-fuzzy metric space if it satisfies the condition (FD).*

As $M(x, y, t)$ is the degree of truth value of the fact that the distance between $x$ and $y$ is less than $t$, for a fixed $0 < \epsilon < 1$, the truth value of the statement "*the actual distance between $x$ and $y$ is* $\inf\{t : M(x, y, t) > \lambda\}$" could be understood as $\lambda$. Suppose the limit of $\inf\{t : M(x, y, t) > \lambda\}$ exists, as $\lambda$ increases to 1, then the limit can be accepted as the actual distance between $x$ and $y$ induced by the fuzzy metric $M$. Thus, we define the following.

**Definition 3.9** (Actual metric induced by a fuzzy metric). *Let* $(X, M)$ *be a FD-fuzzy metric space. We define the actual metric induced by the fuzzy metric M by* $d_M(x, y) = \lim_{\lambda \to 1} \Delta_{M;\lambda}(x, y)$, *provided the limit exists, for all* $x, y \in X$.

**Lemma 3.10.** *Let* $(X, M)$ *be a fuzzy metric space. If* $0 < \lambda_1 < \lambda_2 < 1$, *then* $\Delta_{M;\lambda_1}(x, y) \leq \delta_{M;\lambda_2}(x, y)$, $\forall x, y \in X$.

*Proof.* If $\delta_{M;\lambda_2}(x, y) < \Delta_{M;\lambda_1}(x, y)$, then we choose $t_0 \in \mathbb{R}$ such that $\delta_{M;\lambda_2}(x, y) < t_0 < \Delta_{M;\lambda_1}(x, y)$. Hence, we have $t_0 \notin (-\infty, \delta_{M;\lambda_2}(x, y)) = \{t : M(x, y, t) < \lambda_2\}$ and $t_0 \notin (\Delta_{M;\lambda_1}(x, y), \infty) = \{t : M(x, y, t) > \lambda_1\}$. Therefore, $\lambda_1 \geq M(x, y, t_0) \geq \lambda_2$, which is a contradiction. Thus the lemma follows. $\square$

**Theorem 3.11.** *Let* $(X, M)$ *be a fuzzy metric space. Then the following are equivalent.*

*(1)* $(X, M)$ *is an (FD) fuzzy metric space.*
*(2)* $\lim_{\lambda \to 1} \delta_{M;\lambda}(x, y)$ *exists for all pairs* $(x, y)$.
*(3)* $\lim_{\lambda \to 1} \Delta_{M;\lambda}(x, y)$ *exists for all pairs* $(x, y)$.

*Proof.* (1) $\Rightarrow$ (2) We recall that Condition (*FD*) states that for every pair $(x, y)$ of points in $X$, there exists $t_{x,y}$ such that $M(x, y, t_{x,y}) = 1 > \lambda$, $\forall \lambda \in (0, 1)$, and hence $\Delta_{M;\lambda}(x, y) = \inf\{t : M(x, y, t) > \lambda\} \leq t_{x,y}$. Since $\delta_{M;\lambda}(x, y) \leq \Delta_{M;\lambda}(x, y) \leq t_{x,y}$, for all $\lambda \in (0, 1)$ and $\delta_{M;\lambda}(x, y)$ increases with $\lambda$ (See Lemma 3.4), we see that $\lim_{\lambda \to 1} \delta_{M;\lambda}(x, y)$ exists.

(2) $\Rightarrow$ (3) Let $\lambda \in (0, 1)$ be arbitrary. We choose $\lambda'$ between $\lambda$ and 1. Then, by Lemma 3.10, we have $\Delta_{M;\lambda}(x, y) \leq \delta_{M;\lambda'}(x, y) \leq \lim_{\lambda \to 1} \delta_{M;\lambda}(x, y)$, as $\delta_{M;\lambda}(x, y)$ increases with $\lambda$. Since $\Delta_{M;\lambda}(x, y)$ increases as $\lambda$ increases and bounded above, we conclude that $\lim_{\lambda \to 1} \Delta_{M;\lambda}(x, y)$ exists.

(3) $\Rightarrow$ (1) Assume that $\lim_{\lambda \to 1} \Delta_{M;\lambda}(x, y)$ exists and let it be $t_0$. Since $\Delta_{M;\lambda}(x, y)$ increases as $\lambda$ increases, we have $t_0 + 1 > \Delta_{M;\lambda}(x, y)$ for all $0 < \lambda < 1$. Hence $M(x, y, t_0 + 1) > \lambda$, for all $0 < \lambda < 1$. Thus, $M(x, y, t_0 + 1) = 1$.

Hence the proof is complete. $\square$

**Corollary 3.12.** *Let* $(X, M)$ *be an (FD)-fuzzy metric space. Then, for every* $x, y \in X$, $\lim_{\lambda \to 1} \delta_{M;\lambda}(x, y) = \lim_{\lambda \to 1} \Delta_{M;\lambda}(x, y)$, $\forall x, y \in X$.

*Proof.* Let $x, y \in X$ be fixed arbitrarily. By the condition (FD), both of the limits exist. Using $\delta_{M;\lambda}(x, y) \leq \Delta_{M;\lambda}(x, y), \forall \lambda \in (0, 1)$, we get that $\lim_{\lambda \to 1} \delta_{M;\lambda}(x, y) \leq \lim_{\lambda \to 1} \Delta_{M;\lambda}(x, y)$. On the other hand, from Lemma 3.10, we have that $\Delta_{M;\lambda_1}(x, y) \leq \delta_{M;\lambda_2}(x, y)$, whenever $\lambda_1 <$ $\lambda_2$. If we allow $\lambda_1$ to 1, then $\lambda_2$ also tends to 1, and hence we get the reverse inequality $\lim_{\lambda_1 \to 1} \Delta_{M;\lambda_1}(x,y) \leq \lim_{\lambda_2 \to 1} \delta_{M;\lambda_2}(x,y)$. □

Thus, all the questions we posed in this section, are answered affirmatively.

## 4. Fuzzification of a crisp metric

If $(X, d)$ is a metric space, then the induced fuzzy metrics $M : X \times X \to \mathbb{R}$ are defined in [6] by

$$M_{m,n,k}(x, y, t) = \frac{kt^n}{kt^n + md(x,y)}, \text{ where } m, n, k \in \mathbb{R}^+.$$

We first prove that the topologies constructed from the fuzzy metrics $M_{m,n,k}$ are same and they are equal to the topology compatible with the given metric $d$. As there are infinitely many such fuzzy metrics, we prefer to find a unique fuzzy metric obtained from the given metric. Since the traditional way of fuzzifying a crisp object is identifying the crisp object as a function with the range $\{0, 1\}$, we now introduce a fuzzy metric with range $\{0, 1\}$ representing the given metric in Theorem 4.2.

**Lemma 4.1.** *For each $m, n, k \in \mathbb{R}^+$, the collection of all open balls with respect to the fuzzy metric $M_{m,n,k}$ is same as that of the original metric $d$.*

*Proof.* Let $(X, d)$ be a metric space. Let $B_{m,n,k}(x, r, \epsilon)$ be the open ball with respect to the fuzzy metric $M_{m,n,k}$ with centre $x$, radius $r$ and fuzzy parameter $\epsilon$, as in Definition 2.7. Then

$$\begin{aligned} B_{m,n,k}(x, r, \epsilon) &= \left\{ y : \frac{kr^n}{kr^n + md(x,y)} > 1 - \epsilon \right\} \\ &= \left\{ y : d(x, y) < \frac{\epsilon}{1-\epsilon} \frac{kr^n}{m} \right\} \\ &= B_d(x, s), \end{aligned}$$

where $s = \frac{\epsilon}{1-\epsilon} \frac{kr^n}{m}$. Conversely, if $B_d(x, \eta)$ is a ball in $(X, d)$, then, for any $m, n, k, \epsilon$, if $r$ is chosen as $r = \left(\frac{1-\epsilon}{\epsilon} \frac{\eta m}{k}\right)^{\frac{1}{n}}$, then we get that $B_d(x, \eta) = B_{m,n,k}(x, r, \epsilon)$. Thus the collection of all open balls with respect to the metric $d$ and the collection of all open balls given by the fuzzy metric $M_{m,n,k}$, are the same, for all $m, n, k$. □

In the following theorem, we introduce a fuzzy metric representing the given metric. Incidentally, it is found similar to the fuzzy metirc constructed from a generalized metric in [14, Example 3.8].

**Theorem 4.2.** *If $(X, d)$ is a metric space, then $M_d : X^2 \times \mathbb{R} \to \{0, 1\}$ defined by $M_d(x, y, t) = \begin{cases} 1 & \text{if } d(x, y) < t \\ 0 & \text{if } d(x, y) \geq t \end{cases}$ is a fuzzy metric on $X$.*

*Proof.* It is clear that
- $M_d(x, y, t) = 1$ for all $t > 0$ if and only if $x = y$.
- $M_d(x, y, t) = 0$ for all $x, y$ and for all $t \leq 0$.
- $M_d(x, y, t) = M_d(y, x, t)$ for all $x, y$.

Let us now prove that

(2) $$\min\{M_d(x, y, t), M_d(y, z, s)\} \leq M_d(x, z, t + s).$$

Let $x, y, z \in X$ be arbitrary. If at least one of $M(x, y, t)$ and $M(y, z, s)$ is 0, then clearly (2) follows. If both of them are 1, then we have $d(x, y) < t$ and $d(y, z) < s$. Since $d$ is a metric we have $d(x,z) \le d(x,y) + d(y,z) < t+s$. Thus by definition $M_d(x,z,t+s) = 1$ and hence (2) follows.

Let $x, y \in X$ be such that $x \ne y$. Then $d(x,y) > 0$. If $c = d(x,y)$, then $M_d(x,y,t) = 0$ on $(-\infty, c]$ and $M_d(x,y,t) = 1$ on $(c,\infty)$. Thus $M_d(x,y,t)$ is left continuous as a function of $t$ at all points and right continuous at 0. □

Hereafter, we mean the fuzzy metric representing the metric $d$ by $M_d$.

**Lemma 4.3.** *For a given $r > 0$ and $\epsilon > 0$, the open ball $B(x,r,\epsilon)$ with respect to $M_d$ is same as the open ball $B_d(x,r)$ with respect to the original metric $d$. Hence, the topologies induced by $d$ and $M$ are same.*

*Proof.* let us find the ball $B(x,r,\epsilon)$ with respect to this fuzzy metric $M_d$.

$$\begin{aligned} B(x,r,\epsilon) &= \{y : M_d(x,y,r) > 1 - \epsilon\} \\ &= \{y : M_d(x,y,r) = 1\} \\ &= \{y : d(x,y) < r\} \\ &= B_d(x,r). \end{aligned}$$

Therefore, the collection of all open balls with respect to $M$ and $d$ are the same, and hence they induce a same topology on $X$. □

From Lemma 4.1 and the above lemma, for any $m, n, k \in \mathbb{R}^+$, we conclude that the topologies induced by $M_{m,n,k}$ and the fuzzy metric with range $\{0, 1\}$ representing $d$ are the same.

**Theorem 4.4** (Consistency between fuzzification and crispification). *If a fuzzy metric $M_d$ is representing a metric $d$, then the upper $\lambda$-metric and lower $\lambda$-metric $\delta_{M;\lambda}$ induced by $M_d$ is same as $d$, for all $\lambda$. Hence, $d = d_{(M_d)}$.*

*Proof.* Let $x, y \in X$ and $\lambda \in (0,1)$. Then, we have

$$\begin{aligned} \Delta_{M;\lambda}(x,y) &= \inf\{t : M(x,y,t) > \lambda\} \\ &= \inf\{t : M(x,y,t) = 1\} \\ &= \inf\{t : t > d(x,y)\} = d(x,y). \end{aligned}$$

and

$$\begin{aligned} \delta_{M;\lambda}(x,y) &= \sup\{t : M(x,y,t) < \lambda\} \\ &= \sup\{t : M(x,y,t) = 0\} \\ &= \sup\{t : t \le d(x,y)\} = d(x,y), \end{aligned}$$

This completes the proof. □

However, if $(X, M)$ is an $FD$-fuzzy metric space, $d_M$ is the metric defined as in Definition 3.9, and $M_{(d_M)}$ is the fuzzy metric defined from $d_M$ as in Theorem 4.2, then unfortunately $M_{(d_M)}$ need not be the original $M$. This statement is justified by the following counterexample.

**Example 4.5.** *If $M : \mathbb{R} \times \mathbb{R} \times \mathbb{R} \to [0,1]$ is defined by*

$$M(x,y,t) = \begin{cases} 0 & \text{if} \quad t \le 0 \\ 1 & \text{if} \quad x = y \text{ and } t > 0 \\ t & \text{if} \quad x \ne y \text{ and } 0 < t \le 1 \\ 1 & \text{if} \quad x \ne y \text{ and } t > 1, \end{cases}$$

*then $M$ is a fuzzy metric on $X$ satisfying the Condition (FD). But as $M_{(d_M)}$ has range $\{0,1\}$, it is not equal to $M$.*

In addition to this, the topologies compatible with $M$ and $M_{(d_M)}$ also need not be the same.

**Example 4.6.** *For $x, y \in \mathbb{R}$ and $t \in \mathbb{R}$, we define*

$$M(x,y,t) = \begin{cases} 0 & \text{if } -\infty < t \leq 0 \\ 1 & \text{if } x = y \text{ and } t > 0 \\ \frac{t}{t+|x-y|} & \text{if } x \neq y \text{ and } 0 < t \leq 2 \\ 1 & \text{if } x \neq y \text{ and } 2 < t. \end{cases}$$

*Then, $(\mathbb{R}, M)$ is an FD-fuzzy metric.*

If $x \neq y$ and $0 < \epsilon < \frac{|x-y|}{2+|x-y|}$, then we can easily obtain that $M(x,y,t) = \frac{t}{t+|x-y|} \leq 1-\epsilon$ for all $t \leq 2$. Therefore, it follows that $d_M(x,y) \geq t$ for all $t \leq 2$, and hence $d_M(x,y) \geq 2$. But $M(x,y,2) = 1 > 1-\epsilon$, which implies that $d_M(x,y) \leq 2$. Thus $d_M(x,y) = 2$ if $x \neq y$ and therefore

$$M_{(d_M)}(x,y,t) = \begin{cases} 0 & \text{if } -\infty < t \leq 0 \\ 1 & \text{if } x = y \text{ and } t > 0 \\ 0 & \text{if } x \neq y \text{ and } 0 < t \leq 2 \\ 1 & \text{if } x \neq y \text{ and } 2 < t \end{cases}$$

Let $B(x_0, r_0, \epsilon_0)$ and $B'(x_0, r_0, \epsilon_0)$ be the open balls with respect to $M$ and $M_{(d_M)}$, respectively. If $r_0 \leq 2$, then $B(x_0, r_0, \epsilon_0) = \left\{x : |x - x_0| < \frac{r_0 \epsilon_0}{1-\epsilon_0}\right\}$, $B'(x_0, r_0, \epsilon_0) = \{x_0\}$ and if $r_0 > 2$, then $B'(x_0, r_0, \epsilon_0) = \mathbb{R} = B'(x_0, r_0, \epsilon_0)$. Thus, the topology of $(\mathbb{R}, M_{(d_M)})$ is strictly finer than that of $(\mathbb{R}, M)$.

Finally, we compare the topologies induced by a fuzzy metric $M$ and that of the metric $d_M$ induced by the fuzzy metric.

**Theorem 4.7.** *Let $(X, M)$ be an FD-metric space. Then, open ball $B(x_0, r_0)$ with respect to the actual metric induced by $d_M$ is contained in $B(x_0, r_0, \epsilon_0)$ with respect to $M$, for all $\epsilon_0 \in (0,1)$. Thus, the topology compatible with $d_M$ is finer than that of $M$.*

*Proof.* Let $\epsilon_0 \in (0,1)$ be arbitrary. If $y \in B(x_0, r_0)$, then $d_M(x_0, y) < r_0$ and hence $\liminf_{\epsilon \to 0}\{t : M(x_0, y, t) > 1-\epsilon\} < r_0$. This implies that $\inf\{t : M(x_0, y, t) > 1-\epsilon_1\} < r_0$, for some $0 < \epsilon_1 < \epsilon_0$. Hence $M(x_0, y, t_0) > 1-\epsilon_1$, for some $t_0 < r_0$. Thus, we get $M(x_0, y, r_0) \geq M(x_0, y, t_0) \geq 1-\epsilon_1 > 1-\epsilon_0$. Hence, $B(x_0, r_0) \subseteq B(x_0, r_0, \epsilon_0)$. This shows that any open set in $(X, M)$ is open in the metric space $(X, d_M)$. □

## 5. Conclusion

In this paper, we have successfully introduced an actual (crisp) metric $d_M$ induced by a fuzzy metric $M$ and proved that the existence of the actual metric induced by a fuzzy metric $M$ is characterized by the $FD$ condition on $M$. We have also provided two different approximations of the metric induced from the fuzzy metric $M$, through upper $\lambda$-metrics $\Delta_{M;\lambda}$ and lower $\lambda$-metrics $\delta_{M;\lambda}$.

On the other hand, considering the fuzzy metrics obtained from a crisp metric, it is natural to identify a metric with a fuzzy metric with range $\{0,1\}$. In this view, we have introduced a fuzzy metric representing a crisp metric and proved that the fuzzy metric is inducing the same topology on $X$ as induced from the original metric. We have also proved that fuzzification of a metric is consistent with the crispification of a fuzzy metric.

Finally, we would like to compare the present work with some of the existing works in the literature. First of all, characterizing crispification of fuzzy metric is stronger than proving the topology induced by an arbitrary fuzzy metric is metrizable [8]. In [16, 17, 1], the metric representing a fuzzy metric $M$ is given, respectively, by

$$
\begin{align}
(3) \quad & d_R(x,y) = \sup\{t \geq 0 : M(x,y,t) \leq 1-t\}, \, \forall x, y \in X, \\
(4) \quad & d_{R_\mu}(x,y) = \sup\{t \geq 0 : M(x,y,t) \leq 1-\mu(t)\}, \, \forall x, y \in X, \\
(5) \quad & d_{R_\alpha}(x,y) = \sup\{t \geq 0 : M(x,y,t) \leq 1-\alpha(t)\}, \, \forall x, y \in X;
\end{align}
$$

where

- $(X, M, *)$ is a fuzzy metric space with the $t$-norm $* \leq *_L$ the Lukasiewicz $t$-norm,
- $(X, M, *)$ is a fuzzy metric space with a $t$-norm and $\mu : \mathbb{R}^+ \to \mathbb{R}^+$ is a continuous function satisfying the properties (i) $\mu(t) = 0 \Leftrightarrow t = 0$, (ii) $\mu(t+s) \geq \mu(t) + \mu(s)$, (iii) $M(x,y,t) > 1 - \mu(t)$, $M(y,z,s) > 1 - \mu(s) \Rightarrow M(x,z,t+s) > 1 - \mu(t+s)$.
- $(X, M, *)$ is a fuzzy metric space and $\alpha : \mathbb{R}^+ \to \mathbb{R}^+$ is an increasing function satisfying the properties (i) $0 < \alpha(t) \leq t$ for all $t \in (0,1)$ and $\alpha(t) > 1$ for all $t > 1$. (ii) $(1 - \alpha(s)) * (1 - \alpha(t)) \geq 1 - \alpha(s+t)$, $\forall s, t \in [0,1]$.

In every construction given above, both the third coordinate $t$ (the distance parameter) of the input of $M$ and the value $M(x,y,t)$ (the truth value of the statement $d(x,y) < t$) of $M$ are dependent, whereas we are considering the general case, where they are independent. Moreover, $d(x,y)$ and $d_{R_\alpha}(x,y)$ are always bounded by 1. Such a restriction is not imposed in our construction.

SOME REMARKS ON METRICS INDUCED BY A FUZZY METRIC 11

R. ROOPKUMAR, DEPARTMENT OF MATHEMATICS, CENTRAL UNIVERSITY OF TAMIL NADU, THIRUVARUR - 610005, INDIA.
*E-mail address*: `roopkumarr@rediffmail.com`

R. VEMBU, DEPARTMENT OF MATHEMATICS, SBK COLLEGE, ARUPPUKOTTAI - 626101, INDIA.
*E-mail address*: `msrvembu@yahoo.co.in`